\documentclass[12pt,leqno]{article}
\usepackage{amsmath,amssymb,amsthm,amscd,dsfont}
\numberwithin{equation}{section} \allowdisplaybreaks
\begin{document}
\newtheorem{theorem}{Theorem}[section]
\newtheorem{defin}{Definition}[section]
\newtheorem{prop}{Proposition}[section]
\newtheorem{corol}{Corollary}[section]
\newtheorem{lemma}{Lemma}[section]
\newtheorem{rem}{Remark}[section]
\newtheorem{example}{Example}[section]
\title{Conformal changes of generalized complex structures}
\author{{\small by}\vspace{2mm}\\Izu Vaisman}
\date{}
\maketitle
{\def\thefootnote{*}\footnotetext[1]%
{{\it 2000 Mathematics Subject Classification: 53C15} .
}}
\begin{center} \begin{minipage}{12cm}
\centerline{\footnotesize Dedicated to Acad. Prof.
Radu Miron on the occasion of his 80-eth birthday}\vspace*{1cm}
A{\footnotesize BSTRACT. A conformal change of $TM\oplus T^*M$ is
a morphism of the form $(X,\alpha)\mapsto(X,e^\tau\alpha)$ $(X\in
TM,\alpha\in T^*M,\tau\in C^\infty(M))$. We characterize the
generalized almost complex and almost Hermitian structures that
are locally conformal to integrable and to generalized K\"ahler
structures, respectively, and give examples of such structures.}
\end{minipage}
\end{center}
\vspace*{5mm}
\section{Introduction} In the last few years, the generalized complex and
K\"ahler structures became an important subject of theoretical
quantum field theory, where they provide new sigma models (e.g.,
\cite{Z}) and allow to express certain supersymmetries (e.g.,
\cite{LRUZ}). This also led to an extensive, purely mathematical
research of the subject (e.g., \cite{Galt}). In this note we discuss
a mathematical question, that of characterizing generalized almost
complex and almost Hermitian structures which become integrable,
respectively, K\"ahler after local conformal changes. The
corresponding classical cases of locally conformal symplectic and
locally conformal K\"ahler structures were studied intensively
(e,g., \cite{DO}). Like in the classical case, the characterization
includes a closed $1$-form $\varpi$ called the Lee form, which
defines the local conformal changes. We construct locally conformal
generalized complex structures and locally conformal generalized
K\"ahler structures, which are not globally conformal, on the Hopf
manifolds and on a product $M\times S^1$ where $M$ is a generalized
Sasakian manifold \cite{V-CRF}. Finally, we discuss the induced
structure on hypersurfaces where the pullback of $\varpi$ vanishes.
\section{Preliminaries}
Throughout the  paper we use the following notation: $M$ is an
$m$-dimensional manifold, $X,Y,..$ are either contravariant vectors
or vector fields, $\alpha,\beta,...$ are either covariant vectors or
$1$-forms, $ \mathcal{X},\mathcal{Y},...$ are pairs
$(X,\alpha),(Y,\beta),...$, $\chi^k(M)$ is the space of $k$-vector
fields, $\Omega^k(M)$ is the space of differential $k$-forms,
$\Gamma$ are spaces of global cross sections of vector bundles, $d$
is the exterior differential and $L$ is the Lie derivative. All the
manifolds and mappings are assumed smooth.

Generalized geometric structures in the sense of Hitchin
\cite{H} are similar to classical structures but defined on the
{\it big tangent bundle} $T^{big}M=TM\oplus T^*M$ with the neutral
metric
$$g((X,\alpha),(Y,\beta))=\frac{1}{2}(\alpha(Y)+\beta(X)) $$ and the
Courant bracket \cite{C}
$$ [(X,\alpha),(Y,\beta)] = ([X,Y],
L_X\beta-L_Y\alpha+\frac{1}{2}d(\alpha(Y)-\beta(X)).$$

A maximal $g$-isotropic subbundle $E$ of $T^{big}M$ (or of the
complexification
$T_c^{big}M=T^{big}M\otimes_{\mathds{R}}\mathds{C}$) is an {\it
almost Dirac structure} and if $\Gamma E$ is closed by the Courant
bracket $E$ is a Dirac structure.

A generalized almost complex structure is a vector bundle
endomorphism $\Phi\in End(T^{big}M)$ that satisfies the following
conditions
$$\Phi^2=-Id,
\;g(\mathcal{X},\Phi\mathcal{Y}) +
g(\Phi\mathcal{X},\mathcal{Y})=0.$$ Furthermore, if the {\it
Nijenhuis torsion} of $\Phi$ vanishes, i.e.,
\begin{equation}\label{Nij}
\mathcal{N}_\Phi(\mathcal{X},\mathcal{Y})=
[\Phi\mathcal{X},\Phi\mathcal{Y}]
-\Phi[\mathcal{X},\Phi\mathcal{Y}] -
\Phi[\Phi\mathcal{X},\mathcal{Y}] +
\Phi^2[\mathcal{X},\mathcal{Y}]=0,\end{equation} where the
brackets are Courant brackets, $\Phi$ is an {\it integrable} or a
{\it generalized complex structure}.

The generalized, almost complex structure $\Phi$ is equivalent
with the pair $(E,\bar E)$ of its $(\pm\sqrt{-1})$-eigenbundles
(the bar denotes complex conjugation), which are complex almost
Dirac structures such that $E\cap\bar E=0$, hence, $\Phi$ may be
defined by $E$. $\Phi$ is integrable iff $E$ is Dirac.

The structure $\Phi$ has the following representation by {\it
classical tensor fields} $$ \Phi\left(
\begin{array}{c}X\vspace{2mm}\\ \alpha \end{array}
\right) = \left(\begin{array}{cc} A&\sharp_\pi\vspace{2mm}\\
\flat_\sigma&-\vspace{-1pt}^{t}A
\end{array}\right) \left( \begin{array}{c}X\vspace{2mm}\\
\alpha \end{array}\right), $$ where
$$\pi\in\chi^2(M),\,\sigma\in\Omega^2(M),\,A\in
End(TM),\,\sharp_\pi\alpha=i(\alpha)\pi,\,\flat_\sigma
X=i(X)\sigma,$$ $t$ denotes transposition, and the following
conditions hold
\begin{equation}\label{alg} \pi(\alpha\circ
A,\beta)=\pi(\alpha,\beta\circ
A),\,\sigma(AX,Y)=\sigma(X,AY),\,A^2=-Id-\sharp_\pi\flat_\sigma.
\end{equation}

In this classical representation the integrability conditions of
$\Phi$ are \cite{{Cr},{V-DGA}}:\vspace{1mm}\\ i) the bivector
field $\pi$ defines a Poisson structure on $M$, i.e.,
$[\pi,\pi]=0$ where the bracket is the Schouten-Nijenhuis bracket
with the sign convention of \cite{V-carte};
\vspace{1mm}\\
ii) the {\it Schouten concomitant} of the pair $(\pi,A)$ vanishes,
i.e., $$R_{(\pi,A)}(\alpha,X) =
\sharp_\pi(L_X(\alpha\circ A)-L_{AX}\alpha)-
(L_{\sharp_\pi\alpha}A)(X)=0;$$ iii) the Nijenhuis tensor of $A$
(defined by (\ref{Nij}) with Lie brackets) satisfies the condition
$$ \mathcal{N}_A(X,Y)  =\sharp_\pi[i(X\wedge
Y)d\sigma];$$ iv) the {\it associated $2$-form}
$\sigma_A(X,Y)=\sigma(AX,Y)$ satisfies the condition
$$
d\sigma_A(X,Y,Z)=\sum_{Cycl(X,Y,Z)}d\sigma(AX,Y,Z).$$

A {\it generalized, Riemannian metric} is a Euclidean (positive
definite) metric $G$ on the bundle $T^{big}M$, which is compatible
with the neutral metric $g$ of $T^{big}M$ in the sense that the
musical isomorphism $$ \sharp_G:T^{big}M=TM\oplus T^*M \rightarrow
T^*M\oplus TM\approx T^{big}M,$$ where $\approx$ is defined by
$(\alpha,X)\leftrightarrows(X,\alpha)$ and
$$
2g(\sharp_G(X,\alpha),(Y,\beta)) = G((X,\alpha),(Y,\beta)),
$$
satisfies the conditions \cite{Galt}
$$
\sharp_G^2=Id,\;
g(\sharp_G(X,\alpha),\sharp_G(Y,\beta)) = g((X,\alpha),(Y,\beta)).
$$

It turns out that a generalized, Riemannian metric is equivalent
with a pair $(\gamma,\psi)$ where $\gamma$ is a classical Riemannian
metric on $M$ and $\psi\in\Omega^2(M)$. More exactly,
$$(\gamma,\psi)\,\leftrightarrows\,
\sharp_G\left(\begin{array}{c}
X\vspace{2mm}\\ \alpha\end{array}\right)=
\left(\begin{array}{cc}\varphi&\sharp_\gamma\vspace{2mm}\\
\flat_\beta &^t\vspace{-1pt}\varphi\end{array}\right)
\left(\begin{array}{c} X\vspace{2mm}\\ \alpha\end{array}\right)
$$ where $\varphi=-\sharp_\gamma\circ\flat_\psi$,
$\flat_\beta=\flat_\gamma\circ(Id-\varphi^2)$  \cite{Galt}.

A {\it generalized almost Hermitian structure} is a pair
$(\Phi,G)$, where $\Phi$ is a generalized almost complex structure
and $G$ is a generalized Riemannian metric, such that the
following skew-symmetry condition holds $$
G(\Phi\mathcal{X},\mathcal{Y})
+G(\mathcal{X},\Phi\mathcal{Y})=0\hspace{3mm}(
\mathcal{X},\mathcal{Y}\in\Gamma T^{big}M).$$
Using the $g$-skew-symmetry of $\Phi$ we see that the previous
condition is equivalent with the commutation condition
$\sharp_G\circ\Phi=\Phi\circ\sharp_G,$ which implies that the pair
$(\Phi^c=\sharp_G\circ\Phi,G)$ ($c$ comes from {\it
complementary}) is a second generalized almost Hermitian structure
that commutes with $\Phi$. A commuting pair $(\Phi,\Phi^c)$
defines $G$ by $\sharp_G=-\Phi\circ\Phi^c$.
\begin{theorem}\label{prGaltK} {\rm\cite{Galt}}
A generalized almost Hermitian structure $(G,\Phi)$ is equivalent
with a quadruple $(\gamma,\psi,J_+,J_-)$, where $\gamma$ is a
classical, Riemannian metric on $M$, $\psi$ is a $2$-form, and
$(\gamma,J_\pm)$ are two classical almost Hermitian structures of
$M$ defined as follows by the matrix of $\Phi$:
\begin{equation}\label{formuleJ}
J_\pm=A+\sharp_\pi\circ\flat_{\psi\pm\gamma}.
\end{equation}
\end{theorem}

The generalized, almost Hermitian manifold $(M,G,\Phi)$ is {\it
generalized, Hermitian} if the structure $\Phi$ is integrable and
{\it generalized, almost K\"ahler } if the complementary structure
$\Phi^c$ is integrable. If both $\Phi$ and $\Phi^c$ are integrable
$(M,G,\Phi)$ is a {\it generalized, K\"ahler} manifold. The
classical structures with the same names yield the simplest
examples.
\begin{theorem}\label{theoremGalt} The structure $(G,\Phi)$ is generalized
K\"ahler iff one of the following hypotheses holds: 1)  $J_{\pm}$
are integrable and one has the equalities
\begin{equation}\label{eqluiGaltieri} d^C_+\omega_+=-d\psi,\,
d^C_-\omega_-=d\psi,\end{equation} where
$\omega_\pm(X,Y)=\gamma(J_\pm X,Y)$ are the K\"ahler forms of the
Hermitian structures $(\gamma,J_\pm)$ and the operators $d^C_\pm$
are defined by the structures $J_\pm$ via the formulas
$$d^C=C^{-1}dC=\sqrt{-1}(\bar{\partial}-\partial)
\;\;(C\lambda=(\sqrt{-1})^{p-q}\lambda,\,\lambda\in\Omega^{p,q}(M));$$
2)  $J_{\pm}$ are integrable and one has the equalities
\begin{equation}\label{genKprinCRF}
(\nabla_XJ_\pm)(Y)=\mp\frac{1}{2}\sharp_\gamma[(i(X\wedge
Y)d\psi)\circ J_\pm+i(X\wedge (J_\pm Y))d\psi],\end{equation} where
$\nabla$ is the Levi-Civita connection of the metric $\gamma$; 3)
the $(3,0)$ and $(0,3)$ type components of $d\psi$ are zero and the
connections
\begin{equation}\label{Bismut}
\nabla^\pm_XY=\nabla_XY\pm\frac{1}{2}\sharp_\gamma[i(X\wedge Y)d\psi]
\end{equation} satisfy the condition $\nabla^\pm J_\pm=0$,
respectively.
\end{theorem} Characterizations 1) and 3) were proven by
Gualtieri \cite{Galt}, where it is also shown that
(\ref{eqluiGaltieri}) is equivalent with
\begin{equation}\label{eqluiGaltieri2}
d\omega_+(J_+X,J_+Y,J_+Z) = -d\omega_-(J_-X,J_-Y,J_-Z)=d\psi(X,Y,Z).
\end{equation} The connections (\ref{Bismut}) are
called the Bismut connections and they are the unique metric
connections with covariant torsion $d\psi$. Characterization 2)
follows by replacing $F_\pm$ by $J_\pm$ in Proposition 4.6 of
\cite{V-CRF}.
\section{Locally conformal integrable structures}
Consider the automorphism $ \mathcal{C}_{\tau}: T^{big}M
\rightarrow T^{big}M$ defined by \cite{{GM},{V-stable},{Wd}}
$$ \mathcal{C}_{\tau} (X,\alpha)= (X, e^\tau\alpha),\hspace{3mm}
\tau\in C^\infty(M).$$
We call it a {\it conformal change} of $T^{big}M$ because it
produces a conformal change of the metric $g$: $$
g(\mathcal{C}_{\tau} (X,\alpha), \mathcal{C}_{\tau} (Y,\beta)) =
e^{\tau} g((X,\alpha),(Y,\beta)). $$ Furthermore, if $\tau$ is
locally constant the change will be called a {\it homothety}.

The natural way to apply a conformal change to any $\Phi\in
End(T^{big}M)$ is by conjugation. In particular, the generalized
almost complex structure $\Phi$ and the generalized Riemannian
metric operator $\sharp_G$ will change as follows
$$
\Phi\mapsto\Phi'=
\mathcal{C}_{-\tau}\circ\Phi\circ\mathcal{C}_\tau,\; \sharp_G\mapsto
\sharp_{G'}=\mathcal{C}_{-\tau}\circ\sharp_G\circ\mathcal{C}_\tau.
$$ Accordingly, one gets
\begin{equation}\label{changeend2}
\left(\begin{array}{cc} A&\sharp_\pi\vspace{2mm}\\
\flat_\sigma&-^t\hspace{-1pt}A\end{array}\right) \mapsto
\left(\begin{array}{cc} A&\sharp_{e^\tau\pi}\vspace{2mm}\\
\flat_{e^{-\tau}\sigma}&-^t\hspace{-1pt}A\end{array}\right),\;
\left(\begin{array}{cc}\varphi&\sharp_\gamma\vspace{2mm}\\
\flat_\beta&^t\hspace{-1pt}\varphi\end{array}\right)\mapsto
\left(\begin{array}{cc}\varphi&\sharp_{e^{-\tau}\gamma}\vspace{2mm}\\
\flat_{e^{-\tau}\beta}&^t\hspace{-1pt}\varphi\end{array}\right)
\end{equation} (the minus sign in $e^{-\tau}\gamma$ is because we
look at $\gamma$ as the covariant tensor of the metric). It
follows that if $G\Leftrightarrow(\gamma,\psi)$ then
$G'\Leftrightarrow(e^{-\tau}\gamma,e^{-\tau}\psi)$.

If $(G,\Phi)$ is a generalized almost Hermitian structure
$(G',\Phi')$ is a generalized almost Hermitian structure too.
Moreover, formulas (\ref{formuleJ}) and (\ref{changeend2}) show
that the corresponding pair of classical Hermitian structures does
not change, i.e., $J'_\pm=J_\pm$.
\begin{defin}\label{integrconform} {\rm A generalized almost complex
structure $\Phi$ is {\it globally conformal integrable} if there
exists a conformal change $ \mathcal{C}_\tau$ such that $\Phi'$ is
integrable. If such changes $ \mathcal{C}_\tau$ exist locally (i.e.,
in a neighborhood of each point), $\Phi$ is a {\it locally conformal
integrable structure}. Similarly, one has notions of {\it(locally)
generalized Hermitian, almost K\"ahler and K\"ahler structures}.}
\end{defin}

We obtain the conditions of conformal integrability by applying
conditions i)-iv) of Section 2 to the tensor fields
$(A,e^\tau\pi,e^{-\tau}\sigma)$. The result is
\begin{prop}\label{condconfint} The generalized almost complex
structure $\Phi$ is globally conformal integrable if there exists
a function $\tau\in C^\infty(M)$ such that $\varpi=d\tau$
satisfies the conditions
\begin{equation}\label{iconform} [\pi,\pi]=-2(\sharp_\pi
\varpi)\wedge\pi,\end{equation}
\begin{equation}\label{iiconform} R_{(\pi,A)}(\alpha,X)=\varpi(AX)
\sharp_\pi\alpha-\varpi(X)A\sharp_\pi\alpha,\end{equation}
\begin{equation}\label{iiiconform} \mathcal{N}_A(X,Y) -
\sharp_\pi[i(X\wedge Y)d\sigma]=
-\sigma(X,Y)\sharp_\pi\varpi\end{equation}
$$+\varpi(X)[(Id+A^2)(Y)]-\varpi(Y)[(Id+A^2)(X)],$$
\begin{equation}\label{ivconform}
d\sigma_A(X,Y,Z)-\sum_{Cycl(X,Y,Z)}d\sigma(AX,Y,Z)
\end{equation}
$$=-[\varpi\wedge\sigma_A+(\varpi\circ A)\wedge\sigma](X,Y,Z).$$
\end{prop} \begin{proof}
Condition i) is $[e^\tau\pi,e^\tau\pi]=0$ and a straightforward
calculation shows its equivalence with (\ref{iconform}). If we use
the formula
$$L_{fX}A=fL_XA+(AX)\otimes df-X\otimes(df\circ A)\;\;(f\in C^\infty(M))$$
in ii) for $\Phi'$ the result is (\ref{iiconform}). Furthermore, a
simple calculation gives the following expression of iii) for
$\Phi'$:
\begin{equation}\label{ptiiiconform} \mathcal{N}_A(X,Y) -
\sharp_\pi[i(X\wedge
Y)d\sigma]=-\sharp_\pi[i(X\wedge Y)(d\tau\wedge
\sigma)]
\end{equation}
$$
=-\sigma(X,Y)\sharp_\pi d\tau
-(X\tau)(\sharp_\pi\circ\flat_\sigma)(Y)
+(Y\tau)(\sharp_\pi\circ\flat_\sigma)(X).$$ In view of (\ref{alg})
this result is equivalent to (\ref{iiiconform}) Finally, the new
associated $2$-form is $e^{-\tau}\sigma_A$ and condition iv) for
$\Phi'$ becomes (\ref{ivconform}).\end{proof}
\begin{prop}\label{homothety} Let $\Phi$ be a generalized complex
structure on $M$ and let $\Phi'$ be obtained by a conformal change
of $\Phi$. Assume that $dim\,M>2$ and that $\Phi$ satisfies one of
the following conditions: 1) $\pi$ is non degenerate, 2) $\forall
x\in M$, $A_x^2\neq-Id$ and $A_x$ has no real eigenvalue, 3)
$rank\,\pi>2$ and $\sigma$ is non degenerate. Then $\Phi'$ is
integrable iff the conformal change is a homothety.\end{prop}
\begin{proof} If $\Phi$ is integrable, $\Phi'$ is
integrable too iff the right hand side of the equalities
(\ref{iconform})-(\ref{ivconform}) vanishes. Condition
$(\sharp_\pi\varpi)\wedge\pi=0$ holds iff either $rank\,\pi=2$ or
$\sharp_\pi\varpi=0$. In case 1) we must have the latter equality,
which also implies $\varpi=0$, and we are done. To discuss case 2),
assume that $d_x\tau\neq0$ and take a vector field $X$ such that
$X\tau\neq0$ on a neighborhood $U_x$. Then, the annulation of the
right hand side of (\ref{iiconform}) yields
$A|_{im\,\sharp_\pi}=fId$ on $U_x$. If we apply this equality to a
$1$-form $\flat_\sigma Y$ where the vector field $Y$ is arbitrary
and use (\ref{alg}), we see that $A|_{U_x}$ satisfies an equation of
the form $$ \mathcal{P}(A)=A^3-fA^2+A-fId=(A-fId)(A^2+Id)=0.$$ Since
$A^2+Id\neq0$, the minimal polynomial of $A$ is either $A-fId$ or $
\mathcal{P}(A)$ and $A$ must have a real eigenvalue. Thus, the
hypothesis of case 2) implies $d\tau=0$ as required. In case 3),
since $rank\,\pi>2$ we have $\sharp_\pi d\tau=0$ and the annulation
of the right hand side of (\ref{iiiconform}) reduces to
\begin{equation}\label{h'3}
(X\tau)\sharp_\pi\flat_\sigma(Y)-(Y\tau)\sharp_\pi\flat_\sigma(X)=0.
\end{equation}
On the other hand, $rank\,\pi>2$ implies that $\forall
X\in\chi^1(M)$ with $\sharp_\pi\flat_\sigma(X)\neq0$ there exists
$\lambda\in\Omega^1(M)$ such that $\sharp_\pi\flat_\sigma
X,\sharp_\pi\lambda$ are linearly independent. If $\sigma$ is non
degenerate we may put $\lambda=\flat_\sigma Y$ and (\ref{h'3})
implies $X\tau=0$. Furthermore, if $\sharp_\pi\flat_\sigma(X)=0$,
(\ref{h'3}) reduces to $(X\tau)\sharp_\pi\flat_\sigma(Y)=0$ for any
$Y$ and we get $X\tau=0$ again. Therefore $d\tau=0$.
\end{proof}

Accordingly, we get the following characterization of the locally
conformal integrable, generalized, almost complex structures.
\begin{theorem}\label{lcintegr} Let $(M,\Phi)$ be a
generalized almost complex manifold that satisfies the hypotheses of
Proposition \ref{homothety}. Then $\Phi$ is locally conformal
integrable iff there exists a closed $1$-form $\varpi\in\Omega^1(M)$
such that conditions {\rm(\ref{iconform})-(\ref{ivconform})} hold.
The structure $\Phi$ is globally conformal integrable iff $\varpi$
is exact.
\end{theorem} \begin{proof} If $\varpi$ exists we have a covering
$M=\cup U_a$ such that $\varpi|_{U_a}=d\tau_a$ for some local
functions $\tau_a$ and $
\mathcal{C}_{-\tau_a}\Phi\mathcal{C}_{\tau_a}$ are integrable.
Conversely, if we have a covering $U_a$ of $M$ with functions
$\tau_a$ such that $ \mathcal{C}_{-\tau_a}\Phi\mathcal{C}_{\tau_a}$
are integrable then Proposition \ref{homothety} shows that
$d\tau_a=d\tau_b$ on $U_\alpha\cap U_\beta$. Thus, the local forms
$d\tau_a$ glue up to the required global closed form $\varpi$. The
last assertion of the theorem is obvious.\end{proof}

Like in the classical case \cite{V1976}, we call $\varpi$ the {\it
Lee form}. It is worth noticing that if $\varpi\wedge d\sigma=0$ the
first equality of (\ref{ptiiiconform}) shows that
$(\pi,A,d\sigma-\varpi\wedge\sigma)$ is a Poisson quasi-Nijenhuis
structure \cite{SX}.

In order to get the characterization of generalized, locally
conformal K\"ahler structures we go from a generalized almost
Hermitian structure $(G,\Phi)$ to the equivalent quadruple
$(\gamma,\psi,J_\pm)$, change it to
$(e^{-\tau}\gamma,e^{-\tau}\psi,J_\pm)$ and ask the latter to
satisfy Gualtieri's conditions (\ref{eqluiGaltieri}). The result
is
\begin{prop}\label{gcK} The generalized almost
Hermitian structure $(G,\Phi)$ is conformal generalized K\"ahler
iff $J_\pm$ are integrable and there exists $\tau\in C^\infty(M)$
such that the form $\varpi=d\tau$ satisfies the conditions
\begin{equation}\label{confKcuG1}
d\psi\pm d^C\omega_\pm= \varpi\wedge\psi\mp(\varpi\circ
J)\wedge\omega_\pm.\end{equation}\end{prop}
\begin{proof} The requirement for $J_\pm$ is clear.
For $\Phi'$, (\ref{eqluiGaltieri2}) becomes
$$
(d\psi-d\tau\wedge\psi)(X,Y,Z)=
\pm(d\omega_\pm-d\tau\wedge\omega_\pm)(J_\pm X,J_\pm Y,J_\pm Z).
$$ If we evaluate the wedge products and take into account that
$\omega_\pm(J_\pm X,J_\pm Y)=\omega_\pm(X,Y)$ we get
(\ref{confKcuG1}).\end{proof}
\begin{prop}\label{homotK} If $M$ is a generalized K\"ahler
manifold of dimension greater than $4$ then a conformal change
leads to a generalized K\"ahler structure iff it is a homothety.
\end{prop} \begin{proof} By (\ref{confKcuG1}) the required
condition is $$ (\varpi)\wedge\psi\mp(\varpi\circ
J)\wedge\omega_\pm=0.$$ This implies $\varpi\wedge(\varpi\circ
J)\wedge\omega_\pm=0$. Since $rank\,\omega_\pm>4$ a well known
Cartan lemma tells that this condition holds iff $\varpi=d\tau=0$.
\end{proof}

As a consequence we get
\begin{theorem}\label{locconfK} If $dim\,M>4$,
the generalized almost Hermitian structure $(\gamma,\psi,J_\pm)$ is
a locally conformal, generalized K\"ahler structure iff $J_\pm$ are
integrable and there exists a closed $1$-form $\varpi$ (the Lee
form) that satisfies condition (\ref{confKcuG1}). The same structure
is globally generalized K\"ahler iff $\varpi$ is exact.\end{theorem}

The proof is the same like for Theorem \ref{lcintegr}.

In order to state some other conditions that are equivalent to
(\ref{confKcuG1}) we recall the Weyl connection defined by a
Riemannian metric $\gamma$ and a closed $1$-form $\varpi$:
$$
\tilde\nabla_XY=\nabla_XY-\frac{1}{2}\varpi(X)Y-
\frac{1}{2}\varpi(Y)X+\frac{1}{2}\gamma(X,Y)\sharp_\gamma\varpi,
$$ where $\nabla$ is the Levi-Civita connection of $\gamma$. The
Weyl connection is the Levi-Civita connection of $e^{-\tau}\gamma$
for the local functions $\tau$ that satisfy $d\tau=\varpi$ and it
is the unique torsionless connection such that
$\tilde\nabla_X\gamma=\varpi(X)\gamma$.
\begin{prop}\label{propWB} In Theorem \ref{locconfK},
condition (\ref{confKcuG1}) may be replaced by each of the following
conditions: i) the Weyl connection satisfies the conditions
\begin{equation}\label{WeylptlocK}
(\tilde\nabla_XJ_\pm)(Y)=\mp\frac{1}{2}\sharp_\gamma\{[i(X\wedge
Y)(d\psi-\varpi\wedge\psi)]\circ J_\pm\end{equation} $$+i[X\wedge
(J_\pm Y)](d\psi-\varpi\wedge\psi)\},$$ ii) the connections
\begin{equation}\label{WBismut}
\tilde\nabla^\pm_XY=\tilde\nabla_XY\pm\frac{1}{2}
\sharp_\gamma[i(X\wedge Y)(d\psi-\varpi\wedge\psi)]
\end{equation} satisfy the condition $\tilde\nabla^\pm J_\pm=0$,
respectively.\end{prop} \begin{proof} Instead of using
(\ref{eqluiGaltieri2}), use (\ref{genKprinCRF}), respectively,
(\ref{Bismut}) to express the fact that
$(e^{-\tau}\gamma,e^{-\tau}\psi,J_\pm)$ is a generalized K\"ahler
structure. The integrability of $J_\pm$ implies the annulation of
the $(3,0),(0,3)$ components of $d\psi$ \cite{Galt}.
\end{proof}

Condition i) of Proposition \ref{propWB} shows that if
$d\psi=\varpi\wedge\psi$ then $(\gamma,J_\pm)$ is a pair of
classical, locally conformal K\"ahler structures with the same
metric and the same Lee form. Condition ii) is interesting
because, at least in the generalized K\"ahler case, it may be
related to physics \cite{GHR}.
\begin{example}\label{varHopf} {\rm
Take $M=\mathds{R}^{2n}\backslash\{0\}\approx
S^{2n-1}\times\mathds{R}$ by the diffeomorphism
$\kappa(x)=(x/||x||,ln||x||/ln\lambda)$
$(x\in\mathds{R}^{2n}\backslash\{0\})$ defined for any choice of
$\lambda\in(0,1)$. Denote by $x^i$ $(i=1,...,2n)$ the natural
coordinates on $\mathds{R}^{2n}$ and consider the symplectic form
$\omega=\sum_{h=1}^ndx^h\wedge dx^{n+h}$ and an arbitrary, constant
$(1,1)$-tensor field $A$ that satisfies the condition
$\omega(AX,Y)=\omega(X,AY)$ (such tensor fields obviously exist).
Then $\omega_A$ (defined like $\sigma_A$) is closed, $(\omega,A)$ is
a Hitchin pair \cite{Cr} and it has a corresponding generalized
complex structure $\Phi$ with the chosen tensor field $A$, the
Poisson bivector field $\pi$ defined by
$\sharp_\pi\circ\flat_\omega=-Id$ and the $2$-form $\sigma$ defined
by $\flat_\sigma=\flat_\omega\circ A^2+\flat_\omega$. If we apply to
$\Phi$ the conformal change $\mathcal{C}_{ln||x||^2}$ we get a
conformal integrable, generalized, almost complex structure $\Phi'$
with the tensor fields $(A,||x||^2\pi,(1/||x||^2)\sigma)$. Now,
consider the quotient $
\mathcal{H}^{2n}=(\mathds{R}^{2n}\backslash\{0\})/\Delta_\lambda$
where $\Delta_\lambda$ is the infinite cyclic group generated by the
transformation $x\mapsto\lambda x$, which is called the Hopf
manifold and where $\kappa$ induces a diffeomorphism
$\mathcal{H}^{2n}\approx S^{2n-1}\times S^1$. It is obvious that
$\Phi'$ is invariant by $\Delta_\lambda$. Hence, there exists an
induced generalized, almost complex structure $\Psi$ on
$\mathcal{H}$ and $\Psi$ is locally conformal integrable via the
conformal changes $\mathcal{C}_{-ln||x||^2}$. The conditions
(\ref{iconform})-(\ref{ivconform}) are satisfied for the closed
$1$-form
$$\varpi=-\frac{2\sum_{i=1}^{2n}x^idx^i}{||x||^2}.$$ Since
$\varpi$ is proportional to the length element of $S^1$ (see the
isomorphism $\kappa$) $\varpi$ is not exact and $\Psi$ is not
globally conformal integrable.}\end{example}
\begin{example}\label{exgconfK} {\rm The Hopf manifold $
\mathcal{H}^{2n}$ $(n>1)$ also has locally conformal generalized
K\"ahler structures that are not globally conformal. Indeed, take
the flat metric $\gamma_0=\sum_{i=1}^{2n}(dx^i)^2$ of
$\mathds{R}^{2n}$, an arbitrary constant $2$-form $\psi_0$, and
two $\gamma_0$-compatible, constant, complex structures $J_\pm$,
for instance
$$J_+(\frac{\partial}{\partial x^h})=
\frac{\partial}{\partial x^{n+h}},\; J_+(\frac{\partial}{\partial x^{n+h}})=
-\frac{\partial}{\partial x^{h}},$$ $$J_-(\frac{\partial}{\partial
x^{2h-1}})=
\frac{\partial}{\partial x^{2h}},\; J_-(\frac{\partial}{\partial x^{2h}})=
-\frac{\partial}{\partial x^{2h-1}},$$ where $h=1,...,n$. The
quadruple $(\gamma_0,\psi_0,J_\pm)$ defines a generalized K\"ahler
structure $(G_0,\Phi)$ on $\mathds{R}^{2n}\backslash\{0\}$ and
$(\gamma_0,\psi_0,J_\pm)\mapsto(\gamma_0/||x||^2,\psi_0/||x||^2,J_\pm)$
produces a conformal generalized K\"ahler structure
$(G'_0,\Phi')$. The latter projects to a locally conformal
generalized K\"ahler structure of $ \mathcal{H}^{2n}$, which
satisfies (\ref{confKcuG1}) for $\varpi=-2dln||x||$, and is not
globally conformal generalized K\"ahler because $\varpi$ is not
exact. If $J_+=J_-$ this example reduces to a well known example
of a classical locally conformal K\"ahler structure that is not
globally conformal K\"ahler
\cite{V1976}. It is also known that the manifold $ \mathcal{H}^4$
has no generalized K\"ahler structures with a constant $J_+$
\cite{Galt}.}\end{example}
\begin{example}\label{genS} {\rm Recall that a generalized Sasakian structure
on a manifold $M$ is equivalent with a pair of classical, normal,
almost contact, metric structures $(F_\pm,Z_\pm,\xi_\pm,\gamma)$
complemented by a pair of forms
$\psi\in\Omega^2(M),\kappa\in\Omega^1(M)$ such that the quadruple
$$(e^{t}(\gamma+dt^2),e^{t}(\psi+\kappa\wedge dt),J_\pm=
F_\pm+dt\otimes Z_\pm-\xi_\pm\otimes\frac{\partial}{\partial t})$$
defines a generalized K\"ahler structure on $M\times\mathds{R}$
\cite{V-CRF}. Thus, the similar quadruple without the factors
$e^t$ defines a conformal generalized K\"ahler structure. The
later is invariant by translations along the factor $\mathds{R}$
and it descends to a locally, not globally, conformal, generalized
K\"ahler structure on $M\times S^1$.}\end{example}

Let $(\gamma,\psi,J_\pm,\varpi)$ be a locally conformal generalized
K\"ahler structure on $M$. A hypersurface $\iota:N\hookrightarrow M$
that satisfies the condition $\iota^*\varpi=0$ will be called a {\it
Lee hypersurface} and we will describe the induced structure of an
orientable, Lee hypersurface.

It is known that any orientable hypersurface of a Hermitian manifold
$(M,\gamma,J)$ has an induced metric, almost contact structure
$(F,Z,\xi)$ such that its $\sqrt{-1}$-eigenbundle is a CR-structure
(e.g., \cite{Bl}). The induced structure is obtained by taking a
normal unit vector field $\nu$ of $N$ and by defining
\begin{equation}\label{indcontact}
Z=-J\nu,\,\xi=\flat_\gamma Z, F|_{TN\cap J(TN)}=J,\, F(Z)=0.
\end{equation}

Following \cite{V-CRF} we can say more about the induced
structure. Indeed, using the normal bundle $\nu N=span\{\nu\}$ it
follows easily that, if we look at $J$ as a generalized complex
structure, the hypersurface $N$ is a CRF-submanifold with the
generalized CRF-structure defined by the classical almost contact
structure (\ref{indcontact}) (see Proposition 2.5 and Definition
3.1 of
\cite{V-CRF}). Thus, the induced structure in not just CR, it is
a classical CRF-structure (see Definition 3.2 of \cite{V-CRF}). From
Proposition 3.1 of \cite{V-CRF} it follows that, in our case, the
supplementary integrability condition (besides the CR condition) is
\begin{equation}\label{condsupl}[Z,F^2X+\sqrt{-1}FX]\in H\oplus Q_c,
\end{equation}
where $H$ is the $\sqrt{-1}$-eigenbundle of $F$ and $Q=span\{Z\}$
is the $0$-eigenbundle of $F$. Using the property $F^3+F=0$ it
follows that (\ref{condsupl}) may be changed to
\begin{equation}\label{condsupl1} F\circ(L_ZF)\circ
F=0,\end{equation} equivalently,
\begin{equation}\label{condsupl2} L_ZF=(\xi\circ L_ZF)\otimes Z.
\end{equation}

Another fact to be noticed is that if $\omega$ is the K\"ahler form
of $(\gamma,J)$ then $\iota^*\omega=\Xi$ where
$\Xi(X,Y)=\gamma(FX,Y)$ $(X,Y\in\chi^1(N))$ is the fundamental form
of the structure $(F,Z,\xi,\gamma)$. In particular, we get
\begin{prop}\label{hyperingK} Any orientable hypersurface of a classical
K\"ahler manifold has an induced classical CRF-structure with a
closed fundamental form.\end{prop}

Back to our subject, the announced result about Lee hypersurfaces
is
\begin{prop}\label{Leehypers} An orientable Lee hypersurface of a
locally conformal, generalized, K\"ahler manifold inherits two
metric, almost contact structures $(F_\pm,Z_\pm,\xi_\pm,\gamma)$
with the fundamental forms $\Xi_\pm$ that satisfy the condition
\begin{equation}\label{Lee1} \varpi(\nu)\Xi_\pm=\pm i(Z)\iota^*(d\psi\pm
d^C\omega_\pm).\end{equation}\end{prop} \begin{proof} Pull back
(\ref{confKcuG1}) by $\iota$, then apply the operator
$i(Z)$.\end{proof}
\hspace*{7.5cm}{\small \begin{tabular}{l} Department of
Mathematics\\ University of Haifa, Israel\\ E-mail:
vaisman@math.haifa.ac.il \end{tabular}}
\end{document}